\documentclass[notitlepage]{article}
\usepackage{amssymb,amsmath,amsthm} \usepackage{hyperref}
\usepackage[latin1]{inputenc}
\usepackage[numbers]{natbib}\usepackage{url}
\newif\ifpdf\ifx\pdfoutput\undefinedles\pdffalse\else\pdfoutput=1\pdftrue\fi\ifpdf\pdfcompresslevel=9\fi
\makeatletter\renewcommand\section{ \@startsection
  {section}{1}{\z@}{-3.5ex \@plus -1ex \@minus -.2ex}{2.3ex
    \@plus.2ex}{\bf\centering\normalsize}}
\renewcommand\subsection{\@startsection{subsection}{2}{\z@}{-3.25ex\@plus
    -1ex \@minus -.2ex}{1.5ex \@plus
    .2ex}{\normalfont\normalsize\bfseries}}
\renewcommand\subsubsection{\@startsection{subsubsection}{3}{\z@}{-3.25ex\@plus
    -1ex \@minus -.2ex}{1.5ex \@plus
    .2ex}{\normalfont\normalsize\bfseries}}
\renewcommand\paragraph{\@startsection{paragraph}{4}{\z@}{3.25ex
    \@plus1ex \@minus.2ex}{-1em}{\normalfont\normalsize\bfseries}}
\renewcommand\subparagraph{\@startsection{subparagraph}{5}{\parindent}{3.25ex
    \@plus1ex \@minus
    .2ex}{-1em}{\normalfont\normalsize\bfseries}}\makeatother
\newcommand{\TheoremText}{Theorem}
\newcommand{\ConditionText}{Condition}
\newcommand{\PropositionText}{Proposition}
\newcommand{\LemmaText}{Lemma} \newcommand{\CorollaryText}{Corollary}
\newcommand{\ProofText}{Proof} \newcommand{\AxiomText}{Axiom}
\newcommand{\DefinitionText}{Definition}
\newcommand{\RemarkText}{Remark} \newcommand{\NoteText}{Note}
\newcommand{\ExampleText}{Example}
\newcommand{\ConventionText}{Convention}
\newcommand{\ExerciseText}{Exercise}
\newcommand{\WarningText}{Warning} \newcommand{\ProblemText}{Problem}

\newtheorem{theorem}{\TheoremText}
\newtheorem{condition}{\ConditionText}

\newtheorem{lemma}{\LemmaText} 

\newtheorem{definition}{\DefinitionText}

\newtheorem{varremark}{\RemarkText} \newtheorem{varnote}{\NoteText}
\newtheorem{varexample}{\ExampleText}
\newtheorem{varconvention}{\ConventionText}
\newtheorem{varwarning}{\WarningText}
\newenvironment{remark}{\begin{varremark}\em}{\em\end{varremark}}

\newenvironment{example}{\begin{varexample}\em}{\em\end{varexample}}

\renewenvironment{proof}{
  \noindent\textbf{\ProofText}\ }{\hspace*{\fill}
  \begin{math}\Box\end{math}\medskip}
\newenvironment{proof*}[1]{
  \noindent\textbf{#1\ }}{\hspace*{\fill}
  \begin{math}\Box\end{math}\medskip}
\newcounter{exercisenr} 

\title{Dynamic State Tameness\thanks{ I want to thank professor M. M.
    Rao for a reading of the first version of this paper and
    suggestions made on it that led to a important improvement.  All remaining errors are mine.  This work was 
    supported by COLCIENCIAS (Colombian national science foundation)
    with grant number 284-2003, and Universidad EAFIT.  }}
\author{Jaime A. Londoño\\
  Departamento de Matemáticas\\
  Universidad Nacional de Colombia, Bogotá, Colombia\\
  jlondono@math.ucr.edu}
\begin{document}
\maketitle
\tolerance=200 \setlength{\emergencystretch}{2em}
\begin{abstract}
  An extension of the idea of state tameness is presented in a dynamic
  framework.  The proposed model for financial markets is rich enough
  to provide analytical tools that are mostly obtained in models that
  arise as the solution of SDEs with deterministic coefficients.  In
  the presented model the augmentation by a shadow stock of the price
  evolution has a Markovian character.  As in a previous paper, the
  results obtained on valuation of European contingent claims and
  American contingent claims do not require the full range of the
  volatility matrix.  Under some additional continuity conditions, the
  conceptual framework provided by the model makes it possible to
  regard the valuation of financial instruments of the European type
  as a particular case of valuation of instruments of American type.
  This provides a unifying framework for the problem of valuation of
  financial instruments.
\end{abstract}

\section[Introduction]{Introduction}\label{sec:intro}
State Tameness was introduced in \citet{Londono04} in a setting of a
general semimartingale process driven by Brownian motions, in order to
give a full characterization of non existence of arbitrage that has an
algebraic appealing character with an economic justification.  It also
provided theorems for valuation of financial instruments of European
and American type.  When analytical tools are needed for the study of
financial markets, as it is the case for the problem of optimal
consumption and investment, the general semimartingale framework is
too weak, and the standard approach is usually to impose very strong
conditions (e.g. deterministic coefficients) in order to obtain a rich
theory (\citet{Karatzas98}).  In this paper we propose a model for
financial markets that captures the characteristic properties that a
general formulation of the problem of optimal consumption and
investment, we believe should have.  We think that the main
contribution of this work is that it provides tools and a framework to
solve dynamic problems as the mentioned above. (See for instance
\citet{Londono2005}.)

The model is inspired by heuristic considerations when the author
tried to formulate the problem of optimal investment and consumption
as explained in \citet{Londono2005}.  In particular, in order to model
a typical consumer that changes preferences for consumption partners,
due for instance to aging, it is needed that the model should allow
for the computation of optimal strategies of investment and
consumption after any given time.  As it is explained in
\citet{Londono2005}, the latter problem reduces to computing the
portfolio that finances a given optimal wealth.  Typically, the key
tool to obtain those optimal portfolios, is the representation theorem
for Brownian martingales as stochastic integrals (see for instance
\citet{Londono04}).  In the framework proposed \citet[Exercise
3.2.10]{hK90} can be used instead.  However, in order to apply the
cited theorem, it is necessary that the randomness underlying the
process after any given time be generated by a Brownian motion that
does not carry any information from the past. It is therefore
unavoidable the use of two-parameter filtrations to model the cited
problem.  Finally, we use the concept of non-arbitrage of state tame
portfolios to restrict the class of models we propose to study.

In few words, the model of Market that we propose is a model where the
underlying source of randomness is a finite dimensional Brownian
motion, where the evolution of the process between any two times
depends on the \emph{evolution} of the Brownian motion and the state
of the process at the initial time of the interval. In this model the
state is characterized by the price of the stocks and the value of a
``shadow stock'' that captures the evolution of the economy (see
Section \ref{sec:model}).  If the model is obtained as the result of
solving some stochastic differential equations of coefficients
satisfying some Lipschitz condition, the proposed model becomes a
Brownian flow of homeomorphisms (see Section \ref{sec:model}).  We
point out that there are important classes of processes that do not
fit in the latter framework.  Processes that do not satisfy the
Lipschitz continuity property, and processes driven by Levy processes
(not necessary continuous), are important examples (see
\citet{Kager_Scheutzow1997}).

Using Arnold's terminology (\citet{Arnold98}), the model of prices
proposed is a crude cocycle over the metric dynamical system defined
by the Brownian motion (see remarks at the end of Section
\ref{sec:definitions}).  Arnold's formulation does not require the
evolution of the process to be adapted to the evolution of the
underlying Brownian motion.  In fact when the random dynamical system
is the solution of an SDE driven by a continuous spatial
$(\mathcal{F}_{s,t})$ forward semimartingale helix, satisfying some
smoothness conditions, both definitions are equivalent (See Section
\ref{sec:definitions}).

Next we briefly describe the contents of the paper.  For the sake of
completeness, in Section \ref{sec:definitions} we recall some spaces
of functions, following \citet{hK90}.  Next, we define the notion of
consistent process and we give some examples.  The author is not aware
of a similar definition in the current literature.  Section
\ref{sec:model} defines the model of prices. In order to provide a
consistent framework the model for the price process of $n$ stocks and
a hidden variable proposed is a $(n+1)$ dimensional process.  Its
$k^{th}$-point motion, for any positive integer $k$, satisfies the
Markov property (see the proof of Theorem 4.2.1 in \citet{hK90}).  In
other words the price of the stocks does not necessarily satisfy a
Markov structure, in our model, but the augmentation by a ``hidden''
variable does.  We also define within the model proposed basic
structures in finance as wealth, income and portfolio structures; an
example is provided.  It should be emphasized that although we try to
be as close as possible in definitions and notation for the above
concepts to \citet{Karatzas98}, all of the above structures need a
precise definition that goes along with our framework, since this is
the first time this model is proposed.  After the conceptual framework
is established the math is rather straightforward.  In Section
\ref{sec:statetameness}, a characterization of non existence of
arbitrage opportunities is given, mainly following the proof given in
\citet{Londono04}.  Sections \ref{sec:europeanvaluation} and
\ref{sec:review-state-american} provide the corresponding theory for
valuation of financial instruments of European and American type.
Although some ideas of the proofs given in \citet{Londono04} can be
easily changed to be adapted to the current model, special care must
be taken with the smoothness in the price variable that the model
implies.  As in the mentioned paper full range of the volatility
matrix is not required.  However the most interesting feature of the
framework proposed within the valuation of derivatives is that under
some additional condition on continuity of the expected values of some
random variables (see Condition \ref{con:continuity}), the problem of
valuation of financial instruments of the European type is a
particular case of the problem of valuation of financial instruments
of the American type (see Theorem \ref{thm:completeness} and Theorem
\ref{thm:american}).
   
\section[Some Definitions]{Some Definitions}\label{sec:definitions}
First we introduce some notation which will be frequently used in this
paper.  Let $\mathbb{D}\subset\mathbb{R}^k$ be a open connected set.
Let $m$ be a non-negative integer.  We denote by
$C^{m,\delta}(\mathbb{D}\colon\mathbb{R}^n)$ the the Fréchet space of
$m$-times continuous differentiable functions whose $m$-order
derivatives are $\delta$-H\"older continuous with semi-norms
$\|f\|_{m,\delta\colon K}$ defined in \citet[Section 3.1]{hK90} where
$K\subset \mathbb{D}$ is a compact set and $0\leq\delta\leq 1$.  In
case $m=0$ (or $\delta=0$) we  denote
$C^{m,\delta}(\mathbb{D}\colon\mathbb{R}^n)$ simply by
$C^{\delta}(\mathbb{D}\colon\mathbb{R}^n)$
($C^m(\mathbb{D}\colon\mathbb{R}^n)$).  We also denote by
$\widetilde{C}^{m,\delta}(\mathbb{D}\colon\mathbb{R}^n)$ the Fréchet
space of continuous functions $g\colon
\mathbb{D}\times\mathbb{D}\to\mathbb{R}^n$ which are $m$-times
continuously differentiable with respect to each variable and whose
$m$-order derivatives with respect to both variables are
$\delta$-H\"older continuous with semi-norms
$\|g\|^{\sim}_{m,\delta\colon K}$, where $K$ is a compact set, as
described in \citet[Section 3.1]{hK90}.  In case $m=0$ we  denote
$\widetilde{C}^m(\mathbb{D}\colon\mathbb{R}^n)$ by
$\widetilde{C}(\mathbb{D}\colon\mathbb{R}^n)$ and
$\widetilde{C}^{m,\delta}(\mathbb{D}\colon\mathbb{R}^n)$ by
$\widetilde{C}^{\delta}(\mathbb{D}\colon\mathbb{R}^n)$.  In case
$\delta=0$ we  use the notations
$\widetilde{C}^{m,\delta}(\mathbb{D}\colon\mathbb{R}^n)$ and
$\widetilde{C}^{m}(\mathbb{D}\colon\mathbb{R}^n)$ interchangeably.

We assume a $d$-dimensional Brownian Motion $\{W(t),\mathcal{F}_t;
0\leq t\leq T\}$ starting at $0$ defined on a complete probability
space $(\Omega,\mathcal{F},\mathbf{P})$ where
$\mathcal{F}=\mathcal{F}_T$ and $\{\mathcal{F}_t, 0\leq t\leq T\}$ is
the $\mathbf{P}$ augmentation by the null sets of the natural
filtration $\mathcal{F}^{W}_t=\sigma(W(s),0\leq s\leq t)$.  Let
$(\mathcal{F}_{s,t})_{s_0}=\left\{\mathcal{F}_{s,t}, s_0\leq s\leq t
  \leq T\right\} $ be the two parameter filtration where
$\mathcal{F}_{s,t}$ is the least sub $\sigma$-field containing all
null sets and $\sigma(W_s(u)\mid s\leq u\leq t)$, where $W_s(u)\equiv
W(u)-W(s)$. In the case that $s_0=0$, we just write
$\left(\mathcal{F}_{s,t}\right)$ as an abbreviation for $\left(\mathcal{F}_{s,t}\right)_{0}$.  See
\citet{Arnold98} and \citet{hK90} for detailed accounts of two
parameter filtrations.

Next, we give some definitions frequently used in this paper.  Let
$0\leq s_0\leq T$ be a fixed number.  We shall say that a family
of processes  $\left\{\varphi(s,t), s_0\leq s\leq t\leq T\right\}$ with values in
some euclidean space is a $\left(\mathcal{F}_{s,t}\right)_{s_0}$
\emph{progressive measurable process with two parameters after time
  $s_0$} if for each $s$ with $s_0\leq s\leq T$, $\left\{\varphi(s,t),
  s\leq t\leq T\right\}$ is a $\mathcal{F}_{s,t}$ progressive
measurable processes.  In addition, if for each $s_0\leq s\leq T$, the
process $\{\varphi(s,t), s\leq t\leq T\}$ is a continuous
$\mathcal{F}_{s,t}$-semimartingale, then we  say that the process
$\{\varphi(s,t), s_0\leq s\leq t\leq T\}$ is a continuous
$(\mathcal{F}_{s,t})_{s_0}$-semimartingale with two parameters.
  
Let $\mathbb{D}\subset\mathbb{R}^k$, be an open set. Let
$\varphi(s,t,x,\omega)$, $s_0\leq s \leq t\leq T$, $x\in\mathbb{D}$ be a
$\mathbb{R}^n$-valued random field on the probability space
$(\Omega,\mathcal{F},\mathbf{P})$.  We  call it a 
\emph{measurable process with two parameters after time $s_0$ with
 values in $\mathbb{R}^n$} if for each $x\in\mathbb{D}$,  
$\varphi(\cdot,\cdot,x)$ is a progressive measurable process
with two parameters after time $s_0$.  We  say
that the  $\varphi$ is a
\emph{$C^{m,\delta}(\mathbb{D}\colon\mathbb{R}^n)$-process with
  two-parameters after time $s_0$} if for each $s_0\leq s\leq T$ 
the process $\varphi_s\colon t\to\varphi(s,t,\cdot)$,  is a measurable
random field with values in
$C^{m,\delta}(\mathbb{D}\colon\mathbb{R}^n)$.  In addition, if
$\varphi(s,t,x)$ is a continuous $(\mathcal{F}_{s,\cdot})$ process
for each $x\in\mathbb{D}$ and $s_0\leq s\leq T$, then we shall say
that $\varphi$ is a \emph{continuous
  $C^{m,\delta}(\mathbb{D}\colon\mathbb{R}^n)$-process with two
  parameters after time $s_0$}.  For a definition of measurable random
fields see \citet{hK90}.
  
Let $\varphi$ be a
$C^{m,\delta}(\mathbb{D}\colon\mathbb{R}^n)$-process with two
parameters after time $s_0$.  Assume that $\varphi(s,t,x)$,
$x\in\mathbb{D}$ is a family of continuous semimartingales decomposed
as $\varphi(s,t,x)=\varphi_{loc}(s,t,x)+\varphi_{fv}(s,t,x)$, where
$\varphi_{loc}$ is a continuous
$C^{m,\delta}(\mathbb{D}\colon\mathbb{R}^n)$-local-martingale with two
parameters, and $\varphi_{fv}$ is a continuous
$C^{m,\delta}(\mathbb{D}\colon\mathbb{R}^n)$-process of bounded
variation with two parameters after time $s_0$.  We  say that
$\varphi$ is a \emph{continuous
  $C^{m,\delta}(\mathbb{D}\colon\mathbb{R}^n)$-semimartingale with two
  parameters after time $s_0$}.
  
For each $0\leq s\leq T$ we define a $\sigma$-field $\mathcal{P}_s$ of
progressive measurable sets after time $s$ as the $\sigma$-field of
sets $P\in\mathcal{B}([s,T])\otimes\mathcal{F}_{s,T}$, the product
$\sigma$-field, such that $\chi_{P}(t,\omega)$, $t\geq s$, is a
$\mathcal{F}_{s,t}$ progressive measurable (in $t$) process, where
$\chi$ is the indicator function. Define the measure $\mu_s$ on
$\mathcal{P}_s$ by $\mu_s(P)=\mathbf{E}\int_s^T\chi_P(s,\omega)\,dt$.
Assume that $\varphi$ is a
$C(\mathbb{D}\colon\mathbb{R}^n)$-semimartingale with two parameters
after time $s_0$, with decomposition
$\varphi=\varphi_{loc}+\varphi_{fv}$ as above.  A pair $(a,b)$ where
$a(s,t,x,y)$ and $b(s,t,x)$ are measurable random
fields $\mathcal{F}_{s,t}$-progressive measurable in $t$, for all
$x,y\in\mathbb{D}$, $s_0\leq s\leq T$, is said to be the \emph{local
  characteristics of $\varphi$}, if $(a(s,\cdot,x,y),b(s,\cdot,x))$
is the local characteristic of
$\varphi_s\equiv\varphi(s,\cdot,\cdot)$ (see \citet{hK90}) for any
$s\leq t\leq T$.  In addition, a pair $(\sigma,b)$ where
$\sigma(s,t,x)$ is a measurable random field with values in
$L(\mathbb{R}^d\colon \mathbb{R}^n)$, where
$L(\mathbb{R}^d\colon\mathbb{R}^n)$ denotes the set of matrices with
size $n\times d$, $(\mathcal{F}_{s,t})$-progressive measurable in $t$,
for all $x\in\mathbb{D}$, $s_0\leq s\leq T$ , and $b$ is as above is
said to be the \emph{volatility and drift processes of $\varphi$} if
  \begin{equation*}
      \varphi_{loc}(s,t,x)(\omega)=\int_s^t\sigma(s,u,x)\,dW_s(u), 
  \end{equation*}
  for all $x$, $s$, $t$ and $\omega$.  Assuming a Brownian filtration,
  if $\varphi$ is a continuous
  $C(\mathbb{D}\colon\mathbb{R}^n)$-semimartingale with two parameters
  after time $s_0$, then there exists a pair $(\sigma,b)$ of
  volatility and drift processes of $\varphi$, as a consequence of
  \citet[Exercise 3.2.10]{hK90}.  It follows that the pair of
  volatility and drift processes of a continuous
  $C(\mathbb{D}\colon\mathbb{R}^n)$-semimartingale with two parameters
  after time $s_0$, is unique in the sense that for each
  $x\in\mathbb{D}$, the processes $\sigma(s,\cdot,x)$, and
  $b(s,\cdot,x)$ are determined uniquely up to $\mu_s$-measure $0$.
  Moreover, if we define
  $a=\sigma\sigma^{\prime}=\left\{\sigma(s,t,x)\sigma(s,t,y) ;
    x,y\in\mathbb{D}, s\leq t\leq T \right\}$, then $\varphi$ is a
  process with local characteristic $(a,b)$, and a similar remark to
  the one made for the uniqueness of the volatility and drift
  processes applies to the uniqueness of the local characteristic.
 
  We  say that \emph{$\varphi$ has local characteristic and drift
    of class ${C}^{m,\delta}(\mathbb{D}\colon\mathbb{R}^n)$} if
  $\varphi_s$ has local characteristic and drift of class
  $\widetilde{C}^{m,\delta}(\mathbb{D}\colon\mathbb{R}^n)$ for any
  time $s\geq s_0$.  Similarly we  say that \emph{$\varphi$ has
    volatility and drift processes of class
    ${C}^{m,\delta}(\mathbb{D}\colon\mathbb{R}^n)$}.  If $\varphi$ is
  a continuous semimartingale with volatility and drift of class
  ${C}^{m,\delta}$, for $\delta>0$, then its local characteristic
  belongs to the class $\widetilde{C}^{m,\delta}$, and it follows by a
  well known result that $\varphi$ has a modification to a continuous
  semimartingale of class $C^{m,\epsilon}$ for any $\epsilon<\delta$.
  (This follows as a consequence of \citet[Theorem 3.1.1]{hK90}).
  Reciprocally, if $\varphi$ is a continuous semimartingale of class
  $C^{m,\epsilon}(\mathbb{D}\colon\mathbb{R}^n)$ it follows (as a
  consequence of \citet[ Exercise 3.2.10 (iii)]{hK90}) that the
  volatility and drift can be chosen to be progressive measurable
  processes of class $C^{m,\delta}$ for any $\delta<\epsilon$.
  
  Let $\varphi(s,t,x), x\in\mathbb{D}$ and $\psi(s,t,x),
  x\in\mathbb{D}$ be measurable processes with two parameters after
  time $s_0$ with values in $\mathbb{R}^n$ and $\mathbb{D}$,
  respectively; in addition, it is assumed that $\psi(s,s,x)=x$ for
  all $x\in\mathbb{D}$, and $0\leq s\leq T$.  We  say that the
  process $\varphi$ is a \emph{$\psi$-consistent process} if for each
  $s_0\leq s\leq s^{\prime}\leq T$ there exists a set
  $N_{s,s^{\prime}}\in \mathcal{P}_{s^{\prime}}$ with
  $\mu_{s^{\prime}}(N_{s,s^{\prime}})=0$, such that
  $\varphi(s,t,x)=\varphi(s^{\prime},t,\psi(s,s^{\prime},x))$ for
  all $(t,\omega)\notin N_{s,s^{\prime}}$ and all $x\in\mathbb{D}$.
  We  say that the process $\varphi$ is a \emph{consistent
    process} if $\varphi$ is a $\varphi$-consistent process.
  
  Let $\tau=\left\{\tau(s,x), x\in\mathbb{D},s_0\leq s\leq T\right\}$
  be a family of stopping times with values in $\left[s_0,T\right]$.
  It is assumed that for each $s_0\leq s\leq T$, $x\in\mathbb{D}$,
  $\tau(s,x)$ is a stopping time relative to the filtration
  $\left\{\mathcal{F}_{s,t}; s\leq t\leq T\right\}$, and that
  $\tau(s,x)(\omega)$ is a measurable random field that  is lower
  semi-continuous with respect to $(s,x)$. We  say that a family $\tau$ as above is a
  \emph{measurable family of stopping times after time $s_0$}; we
   say that the random field $\psi(s,t,x), s_0\leq s\leq t\leq
  \tau(s,x), x\in\mathbb{D}$ is a \emph{measurable process of
    two parameters after time $s_0$ with random time $\tau$}, if
  $\psi_{\tau}=\{\psi(s,\tau(s,x)\wedge t,x), s_0\leq t\leq T ; x,
  s_0\leq s\leq T\}$ where $s\wedge t=\min\{s,t\}$ is a measurable
  process with two parameters after time $s_0$; we  say that the
  family $\tau$ is a \emph{$\psi$-consistent family of stopping times
    after time $s_0$}, if for each $s_0\leq s\leq T$ there exist
  $N_{s}\in\mathcal{P}_s$, $\mu_s(N_s)=0$ with
  $\tau_s(x)=\tau_{t\wedge\tau}(\psi(s,t\wedge\tau,x))$ for all
  $(t,\omega)\notin N_s$ and all $x$, and in this case we  say
  that $(\psi,\tau)$ is a \emph{consistent stopping structure}.  We
  say that the consistent stopping structure $(\psi,\tau)$ is of
  class $C^{m,\delta}(\mathbb{D}\colon\mathbb{R}^k)$ if $\psi_{\tau}$
  is a process of class $C^{m,\delta}(\mathbb{D}\colon\mathbb{R}^k)$.
  Given a consistent stopping structure $(\psi,\tau)$, we  say
  that a family of $\mathbb{R}^n$-valued processes
  $\varphi=\{\varphi(s,t,x), s\leq t\leq \tau_s(x) ; x\in\mathbb{D},
  s_0\leq s\leq T\}$ is a \emph{$\psi$-consistent process with random
    time $\tau$}, if $\varphi_{\tau}$ is a $\psi_{\tau}$-consistent
  measurable process with two parameters after time $s_0$. Similarly,
  we  say that $\varphi$ is a process of class
  $C^{m^{\prime},\delta^{\prime}}(\mathbb{D}\colon\mathbb{R}^n)$ if
  $\varphi_{\tau}$ is a process of the same class.

  Before we start discussing a model for financial markets using the
  above terminology, we do a digression on how this fits in the
  framework of the theory of Random Dynamical Systems introduced by L.
  Arnold and his school (see \citet{Arnold98}).  The notation on the
  following paragraph is local to it.  Let
  $(\Omega,\mathcal{F},\mathbf{P})$ be the probability space defined
  above.  Let $\theta\colon\mathbb{R}\times\Omega\to\Omega$ be the
  $P$-preserving flow on $\Omega$ defined by $\theta(t, \omega)\equiv
  W_{\cdot}-W_t$; without loss of generality we assume that $\theta$
  is defined on $\mathbb{R}\times\Omega$.  Assume that
  $F\colon\mathbb{R}^n\times\mathbb{R}\times\Omega\to\mathbb{R}^d$ is
  a continuous spatial helix forward $(\mathcal{F}_{s,t})$
  semimartingale with forward local characteristic of class
  $\tilde{C}^{m,\delta}$, for $m\geq1$ and $\delta >0$. Assume that
  $\varphi(s,t,\cdot)(\omega)$ is the trajectory random field of the
  differential equation
\[
d\varphi(s,t,x)=F(\varphi(s,t,x),dt),\qquad \varphi(s,s,x)=x.
\]
It is known that $\varphi$ has a jointly measurable modification (also
denoted by $\varphi$), such that $(\varphi,\theta)$ is a perfect
cocycle, with the property that
\[
\varphi(s,s+t)(\omega)=\varphi(0,t)(\theta_s\omega).
\]
It follows that the above structure is both a perfect cocycle and a
$\varphi$ consistent process.  More details can be found in
\citet[Theorem 2.1]{Salah_Eldin_Scheutzow1999}.

We believe that the proposed class of processes has desirable
properties that arise in many areas. We are motivated by the finding
of an appropriate framework for the problem of optimal consumption and
investment in finance.  An elaborated discussion on the heuristics
behind our approach can be found in \citet{Londono2005}.

\section[A mathematical formulation of the model.]{A mathematical formulation of the model.}\label{sec:model}
Fix $0\leq s_0\leq T$; we assume $n+1$ stocks whose \emph{evolution
  price process} $P$ is a consistent $C(\mathbb{R}_+^{n+1}\colon
\mathbb{R}_+^{n+1})$-semimartingale with volatility and drift
processes of class
$C^{0,\epsilon}(\mathbb{R}_+^{n+1}\colon\mathbb{R}^{n+1})$ for some
$\epsilon>0$, where $\mathbb{R}_+$ denotes the set of real positive
numbers. For $0\leq i\leq n$ we define \emph{the price per-share
  process for the $i$-stock}, $P_i$, to be the $P$-consistent
$C(\mathbb{R}_+^{n+1}\colon \mathbb{R}_+)$-semimartingale process
$P_i=\left\{P_i(s,t,p)=\pi_i\circ P(s,t,p), p\in \mathbb{R}^{n+1},
  s_0\leq s\leq t\leq T\right\}$ where $\pi_i$ denotes the projection
on the $i$-component; it follows (see e.g., \citet{Karatzas98}, and
\citet[Exercise 3.2.10]{hK90}) that the evolution price-per-share
process for each stock obeys the differential equations
\begin{eqnarray}\label{eq:iprice_stock}
  dP_i(s,t,p)=P_i(s,t,p)\left[b_i(s,t,p)dt + \sum_{1\leq
        j\leq d}\sigma_{ij}(s,t,p)\,dW^j_s(t)\right]  \\\nonumber
  P_i(s,s,p)=p_{i}, i=1,\ldots,n
\end{eqnarray}
where $W^j_s(t)=W^j(t)-W^j(s)$, and,
\begin{eqnarray}\label{eq:price_shadow_stock}
  dP_0(s,t,p)=P_0(s,t,p)\left[-r(s,t,p)dt - \sum_{1\leq
        j\leq d}\theta_{j}(s,t,p)\,dW^j_s(t)\right]  \\\nonumber
  P_0(s,s,p)=p_0
\end{eqnarray}
for some progressive measurable $P$-consistent processes with two
parameters after time $s_0$ $b_i$, $\sigma_{i,j}$, $r$ and $\theta_i$
of class $C^{0,\delta}(\mathbb{R}_+^{n+1}\colon\mathbb{R}^{n+1})$, for
any $\epsilon>\delta>0$.  We  say that $b_i$ (for $i=1\ldots,n$)
is the \emph{rate of return processes for the $i$-stock}, and that
$\sigma_{i,j}$ is the \emph{(i,j) volatility coefficient processes
  (for $i=1,\ldots,n$, $j=1,\ldots,d$}).

Let us define $F_s(p,t)=(F^0_s(p,t),\cdots,F^n_s(p,t))$ where
\begin{equation}
  \label{eq:Integrands_prices}
F_s^i(p,t)\equiv p_i\int_s^tb_i(s,u,p)\,du + p_i\int_s^t\sum_{1\leq
  j\leq d}\sigma_{ij}(s,u,p)\,dW^j_s(u)
\end{equation}
for $i=1\cdots,n$, and,
\begin{equation}
  \label{eq:Integrands_shadow}
F_s^0(p,t)\equiv -p_0\int_s^tr(s,u,p)du - p_0\int_s^t\sum_{1\leq
  j\leq d}\theta_{j}(s,u,p)\,dW^j_s(u).
\end{equation}
It follows that $P$ is the unique solution of the integral equation
\[
P(s,t,p)=p+\int_s^tF_s(P(s,u,p),du).
\] 

If it is assumed that $b_i(s,t,p)=b_i(p,t)$,
$\sigma_{i,j}(s,t,p)=\sigma_{i,j}(p,t)$, $r(s,t,p)=r(p,t)$, and
$\theta_i(s,t,p)=\theta(p,t)$ where the functions $b_i$,
$\sigma_{i,j}$, $r$, and $\theta$ are deterministic functions for all
$i$, $j$ that are jointly continuous, and Lipschitz continuous in $p$,
then   there is a version
of $P$ that is a forward stochastic flow of homeomorphisms (\citet[Theorem 4.5.1 and Theorem 4.7.1]{hK90}).  We also
observe from \citet[Lemma 4.5.6]{hK90} that Kolmogorov's criterion for
continuous random fields (\citet[Theorem 1.4.1]{hK90}) implies that if
$b$, $\sigma$, $\theta$, and $r$ are of class $C^{0,1}$ (in the price
variable) then $P$ is continuous in $(p,s,t)$.  Indeed, if $b$,
$\sigma$, $\theta$, $r$ are of class $C^{k,\delta}$ for $k\geq 1$ and
$\delta>0$, then \citet[Theorem 4.6.5]{hK90} implies that a version of
$P$ can be chosen that is a forward stochastic flow of
$C^k$-diffeomorphisms.

The meaning of $r(s,\cdot)$, and $\theta_i(s,\cdot)$, for
$i=1,\cdots,d$, $s_0\leq s\leq T$ is as follows; we assume that there
is an imaginary or \emph{shadow stock} whose price evolution is given
by the differential equation \eqref{eq:price_shadow_stock}; additional
conditions are also imposed on $\sigma$, and $r$, as explained below.

We also assume that there is a two parameter \emph{yield process} for
the $i$-th stock, $1\leq i\leq n$, $\left\{Y_i(s,t,p),
  p\in\mathbb{R}_+^{n+1}, s_0\leq s\leq t\leq T\right\}$ that is a continuous $P$-consistent
$C(\mathbb{R}_+^{n+1}\colon\mathbb{R})$-semimartingale.  We assume
that the yield process satisfies the following differential equation
\[
dY_i(s,t,p)=dP_i(s,t,p)+ P_i(s,t,p)\delta_i(s,t,p)\,dt, \qquad
Y_i(s,s,p)=0
\]
for $p\in \mathbb{R}_+^{n+1}$, $s_0\leq s\leq t\leq T$ where
$\delta_i=\left\{\delta_i(s,t,p), s_0\leq s\leq t\leq T, p \in\mathbb{R}_+^{n+1}\right\}$ is a
progressive measurable $P$-consistent process with two parameters
after time $s_0$ of class $C^{0,\epsilon^{\prime}}$, for some
$\epsilon^{\prime}>0$.  We  say that $\delta_i$ is a
\emph{dividend rate processes for the $i$-th stock}.

In addition, it is assumed that the random field  $\theta=\{\theta(s,t,p),
s_0\leq s\leq t\leq T, p\in \mathbb{R}_+^{n+1}\}$ where
$\theta^{\prime}(s,t,p)=\left(\theta_1(s,t),\cdots,\theta_d(s,t)\right)$
for $0\leq s\leq t\leq T$, $p\in \mathbb{R}_+^{n+1}$ is the process
$\theta(s,\cdot,p)\in \ker^{\perp}(\sigma(s,\cdot,p))$, of class
$C^{0,\delta}$ for any $\delta<\min(\epsilon,\epsilon^{\prime})$
(where $\ker^{\perp}(\sigma(s,\cdot,p)$ denotes the orthogonal
complement of the kernel of $\sigma(s,\cdot,p)$) such that
\begin{multline}\label{E:wviability}
  b(s,t,p)+\delta(s,t,p)-r(s,t,p)\mathbf{1}_n\\
  -proj_{\ker(\sigma^{\prime}(s,t,p))}(b(s,t,p)+\delta(s,t,p)-r(s,t,p)\mathbf{1}_n)\\
  =\sigma(s,t,p)\theta(s,t,p)
\end{multline}
a.e. $\mu_s$, for all $p\in\mathbb{R}_+^{n+1}$, and $0\leq s\leq T$,
where $\mathbf{1}_n^{\prime}=(1,\cdots,1)\in\mathbb{R}^n$.  Let us
observe that although $\theta $ is always well defined, in the sense
that there is a progressive measurable process that satisfies the
above equation, it is usually not a process of class $C^{0,\delta}$
for some $\delta$.  Examples of markets where $\theta$ is of class
$C^{0,\delta}$ for some $\delta$, are those for which
$Im(\sigma(s,t,p))$ is a fixed subspace.  We  say that the
process $\theta$ is the \emph{market price of risk}.
  
The $P$-consistent $C(\mathbb{R}_+^{n+1}\colon\mathbb{R}_+)$-process
$B=\left\{B(s,t,p)\right\}$ of bounded variation,
whose evolution $B(s,\cdot,p)$, $p\in\mathbb{R}_+^{n+1}$, $0\leq
s\leq T$ is given by the stochastic differential equation
\begin{equation}\label{eq:Bond}
  dB(s,t,p)=B(s,t,p)r(s,t,p)dt,\qquad B(s,s,p)=1, \text{ for
  }s_0\leq s\leq t\leq T
\end{equation}
will be called the \emph{bond price process} and we  say that $r$
is the \emph{interest rate process}.

We shall say that $\mathcal{M}=(P,b,\sigma,\delta,r,p^0)$ is a
\emph{financial market with terminal time $T$ and initial time $s_0$}
if $b=(b_1,\ldots,b_n)$ is a vector of rate of return processes,
$\sigma=(\sigma_{i,j})$ is a matrix of volatility coefficient
processes, $\delta=(\delta_1,\ldots,\delta_n)$ is vector of dividend
rate processes, $r$ is an interest rate process as explained above,
and $p^0\in\mathbb{R}_+^{n+1}$ is a vector of initial prices.  Let us
observe that if $\mathcal{M}$ is a financial market with initial time
$0$ and terminal time $T$, then for any $0\leq T_0\leq T$ the
restrictions (defined in the obvious way) of $b$, $\sigma$, $\delta$,
and $r$ to the parameter set $T_0\leq s\leq t\leq T$ with respect to
the (two-parameter) filtration $(\mathcal{F}_{s,t})_{T_0}$, along with
any $p\in\mathbb{R}_+^{n+1}$, is a financial market with
initial time $T_0$ and terminal time $T$.  

We define the \emph{state price density process} to be the continuous
$C(\mathbb{R}^{n+1}_+\colon\mathbb{R}_+)$-semimartingale process
defined by
\begin{equation*}
  H(s,t,p)=B^{-1}(s,t,p)Z(s,t,p)\qquad\text{for } p\in
  \mathbb{R}_+^{n+1}, 0\leq s\leq t \leq T
\end{equation*}
where
\begin{equation*}
  Z(s,t,p)=\exp\left\{-\int_{s}^t\theta^{\prime}(s,u,p)\,dW_{s}(u)\ 
-\frac{1}{2}\int_{s}^t\left\|\theta(s,u,p)\right\|^2\,du\right\}
\end{equation*}
for $0\leq s\leq t\leq T$, and $B^{-1}(s,t,p)=1/B(s,t,p)$.  From the
given definitions it follows that the processes $H$ and $P_0$ are
related by the equations
\[
P_0(s,t,p)=p_0H(s,t,p), \qquad\text{for } p\in \mathbb{R}_+^{n+1},
0\leq s\leq t\leq T.
\]

Fix $s_0\in\left[0,T\right]$.  Assume that $\tau=\{\tau_s(x,p);
s_0\leq s\leq T, x\in\mathbb{R}, p\in\mathbb{R}_+^{n+1}\}$ is a
measurable family of stopping times after time $s_0$. \emph{A wealth
  structure after time $s_0$} is a triple $(X, \tau, x_0)$, where
$x_0\in\mathbb{R}$, and $X=\{X(s,t,x,p); x\in\mathbb{R},
p\in\mathbb{R}_+^{n+1}, s\leq t\leq \tau_s(x,p)\}$, is a family of
continuous semimartingale processes with the property that
$((X,P),\tau)$ is a consistent stopping structure of class
$C^{0,\epsilon}(\mathbb{R}\times\mathbb{R}_+^{n+1}\colon\mathbb{R}\times\mathbb{R}_+^{n+1})$
for some $\epsilon=\epsilon^{X}>0$ where $(X,P)$ is the continuous
process with two parameters after time $s_0$, with random time $\tau$
defined as
\[
(X,P)=\left\{(X(s,t,x,p),P(s,t,p)), x\in\mathbb{R}, p\in
  \mathbb{R}_+^{n+1}, s\leq t\leq \tau_s(x,p)\right\}.
\]
We  say that $x_0$ is the \emph{initial value for the wealth
  process}, and we  say that $(X,\tau)$ is a \emph{wealth
  evolution structure}; we shall denote this by writing
$(X,\tau)\in\mathcal{X}(\mathcal{M})$.  Note that the family of
processes $(P_i,T)$, for $i=0,\cdots,n$ and $(\{xB(s,t,p), s_0\leq
s\leq t\leq T,x\in\mathbb{R}, p\in \mathbb{R}_+^{n+1}\},T)$ are wealth evolution structures, as any
reasonable definition should account for.

Next we define a portfolio-income structure. Assume that $(X, \tau)$
is a wealth evolution structure after time $s_0$.  Let $\Gamma$ be a
continuous semimartingale process with random time
$\tau$ after time $s_0$ of class
$C^{0,\epsilon}(\mathbb{R}\times\mathbb{R}_+^{n+1}\colon\mathbb{R})$
(where the positive number $\epsilon$ depends on $\Gamma$) with the
property that $\Gamma(s,s,x,p)=0$  and
\[
\Gamma(s,t^{\prime},x,p)+\Gamma(t^{\prime},t,X(s,t^{\prime},x,p),P(s,t^{\prime},p))=\Gamma(s,t,x,p)
\]
for all $x\in\mathbb{R}$,
$p\in\mathbb{R}_+^{n+1}$, and $s_0\leq s\leq \tau_s(x,p)$. 
We  say that a process $\Gamma$ as above is an \emph{income
  evolution structure for the wealth evolution structure $(X,\tau)$},
and we say that $(X,\Gamma,\tau)$ is a \emph{wealth and income
  evolution structure}.  If $\Gamma(s,t,x,p)\leq 0$ for all $x$, $p$,
$s_0\leq s\leq \tau_s(x,p)$ we  say that $\Gamma$ is a
\emph{consumption evolution structure for the wealth evolution
  structure $(X,\tau)$}.  Let
$(\pi_0,\pi)=\{(\pi_0(s,t,x,p),\pi(s,t,x,p)); x\in\mathbb{R},
p\in\mathbb{R}_+^{n+1}, s_0\leq s\leq t\leq \tau_s(x,p)\}$ be a
$(X,P)$-consistent progressive measurable process of class
$C^{0,\epsilon}$ for some $\epsilon>0$ with random time $\tau$, and
$\pi_o+\pi^{\prime}\mathbf{1}_n=X$ satisfying
\begin{multline}\label{eq:wealthportfolio_process} 
  B^{-1}(s,t,p)X(s,t,x,p)=x+\int_{s}^tB^{-1}(s,u,p)\,d\Gamma(s,u,x,p)\\
  +\int_{s}^tB^{-1}(s,u,p)\pi^{\prime}[x,p](s,u,x,p)\sigma(s,u,p)\,dW_{s}(u) \\
  +\int_s^tB^{-1}(s,u,p)\pi^{\prime}(s,u,x,p)(b(s,u,p)+\delta(s,u,p)-r(s,u,p)\mathbf{1}_n)\,du
\end{multline}
for all $x\in\mathbb{R}$, $s_0\leq s \leq \tau_s(x,p)$,
$p\in\mathbb{R}^{n+1}_+$.  We  say that
$((\pi_0,\pi),\Gamma,\tau)$ as above is a \emph{portfolio evolution
  structure with random time $\tau$ after time $s_0$, financed by the
  income $\Gamma$}.  If $x_0\in\mathbb{R}$, we  say that
$((\pi_0,\pi),\Gamma,\tau,x_0)$ is a \emph{portfolio structure for the
  random time $\tau$ after time $s_0$, financed by the income $\Gamma$
  with initial wealth $x_0$}.  We  say that a \emph{wealth
  evolution structure $(X,\tau)\in\mathcal{X}(\mathcal{M})$ after time
  $s_0$ is financed by the income structure $\Gamma$}, if there exists
a portfolio evolution structure $((\pi_0,\pi),\Gamma,\tau)$ with
random time $\tau$ after time $s_0$ with
$\pi_0+\pi^{\prime}\mathbf{1}_n=X$.  In this case we  say that
$(X,\Gamma,\tau)$ is a \emph{hedgeable wealth-income structure after
  time $s_0$}. Whenever $s_0=0$ in any of the structures above, we
omit explicit reference to the initial time.  In case
$(X,\Gamma,\tau)$ is a hedgeable wealth-income structure after time
$s_0$ and $\Gamma\equiv 0$, we  say that the wealth-income
structure $(X,\tau)$ is a \emph{self-financed wealth evolution
  structure}.

\begin{example}\label{ex:construction_wealth}   
  It is possible to construct structures with the above
  characteristics.  For instance, the simplest example is when the
  ``relevant'' parameters of the model are deterministic functions.
  Assume $b(p,t)$, $\sigma(p,t)$, $\delta(p,t)$, $r(p,t)$, and
  $\theta(p,t)$ are continuous functions in
  $(p,t)\in\mathbb{R}_+^{n+1}\times [0,T]$, that are locally Lipschitz
  continuous in $p$ with values in $\mathbb{R}^{n}$,
  $L(\mathbb{R}^d\colon\mathbb{R}^n)$, $\mathbb{R}^n$, $\mathbb{R}$,
  and $\mathbb{R}^d$, respectively, with
\begin{multline*}
  b(p,t)+\delta(p,t)-r(p,t)\mathbf{1}_n\\
  -proj_{\ker(\sigma^{\prime}(p,t))}(b(p,t)+\delta(p,t)-r(p,t)\mathbf{1}_n)\\
  =\sigma(p,t)\theta(p,t)
\end{multline*}
for all $p$, $t$, where $\theta(p,\cdot)$ belongs to the orthogonal
complement of the kernel of $\sigma(p,\cdot)$.  It follows that the
above functions define a financial market with terminal time $T$.
Assume that $\Gamma$ is the continuous
$C(\mathbb{R}\times\mathbb{R}^{n+1}_+\colon\mathbb{R})$-semimartingale
\[
\Gamma(x,p,s,t)=\int_s^t
b_{\Gamma}(x,p,u)\,du+\int_s^t\sigma_{\Gamma}(x,p,t)\,dW_s(u)
\]
where $b_{\Gamma}$ and $\sigma_{\Gamma}$ are continuous functions in
$(x,p,t)\in\mathbb{R}\times\mathbb{R}^{n+1}\times [0,T]$, which are
locally Lipschitz continuous in $(x,p)$.  Moreover, assume a
$\mathbb{R}^n$-valued function $\pi(x,p,t)$ that is continuous in
$(x,p,t)$ and locally Lipschitz continuous in $(x,p)$.  It follows
that the solution of
\begin{multline*}
  X(s,t,x,p)=x+\int_s^tX(s,t,x,p)r(p,u)\,du\\
  +\int_s^t\left[b_{\Gamma}(x,p,u)\,du+\sigma_{\Gamma}(x,p,u)\,
    dW_{s}(u)\right]\\
  +\int_{s}^t\pi^{\prime}(x,p,u)\left[\sigma(p,u) dW_{s}(u)+ (b(p,u)+
    \delta(p,u)-r(p,u)\mathbf{1}_n)\,du\right]
\end{multline*}
defines a unique solution process $X$ on any random time interval
$[s,\tau^{\prime}]$ before explosion, for each $s\geq s_0$.

More precisely, let
$F_s(x,p,t)\equiv(F_s^0(x,p,t),\cdots,F_s^{n+1}(x,p,t))$ where
\begin{multline*}
  F_s^{n+1}(x,p,t)=\int_s^t(x\, r(p,u)+b_{\Gamma}(x,p,u)+b(p,u)+
  \delta(p,u)-r(p,u)\mathbf{1}_n)\,du\\
  + \int_s^t\left[\sigma_{\Gamma}(x,p,u)
    +\pi^{\prime}(x,p,u)\sigma(p,u)\right] dW_{s}(u)
\end{multline*}
and the processes $F_s^0,\cdots,F_s^{n}$ are defined by equations
\eqref{eq:Integrands_prices} and \eqref{eq:Integrands_shadow}, where
$b_i(s,t,p)=b_i(t,p)$, $\sigma_{i,j}(s,t,p)=\sigma_{i,j}(t,p)$,
$r(s,t,p)=r(t,p)$, and $\theta_j(s,t,p)=\theta_j(t,p)$ for all $i$,
$j$, $s$, and $t$.  Let $((X,P),\tau)$ be the consistent stopping
structure defined as the local solution to the equation
\[
\tilde{X}(s,t,p)=p+\int_s^tF_s(\tilde{X}(s,u,p),du)
\]
where $\tilde{X}=(X,P)$, and $\tau$ is the explosion time.  Let
$(\tau_n)$ be a sequence of $\tilde{X}$-consistent families of
stopping times increasing to $\tau$ with the property that
$\tau_n(s,x,p)<\tau(s,x,p)$ and $\tau_n(s,x,p)\uparrow\tau(s,x,p)$
holds a.s..  It follows that $(X, \Gamma, \tau^{\prime})$ is a
wealth-income structure for any stopping time for which there exists
$n$, with $s_0\leq \tau^{\prime}\leq\tau_n$.  The meaning is natural.
Assume that a market is given and that an investor has an income after
time $s\geq s_0$, that depends on his current wealth $x$ and the state
of the economy (reflected by the value of the stocks
$p_{1,\cdots,n}\equiv (p_1,\cdots,p_n)$ and shadow stock value $p_0$).
Assume also that the increase of the wealth is about
\[
d\Gamma(x,p,s,s+\delta)\approx
b_{\Gamma}(x,p,s)\,d\delta+\sigma_{\Gamma}(x,p,u)\,dW_s(s+\delta)
\]
at any time $s$ after time $s_0$ for small $\delta$, where
$p=(p_0,\cdots,p_n)$.  Moreover assume that the investor has a
strategy to invest in stocks that consists in holding at time $s$,
$\pi(x,p,s)$ worth of stocks (when he has $x$ dollars of wealth, the
vector of prices of stocks is $(p_1,\cdots,p_n)$, and the price of the
shadow stock is $p_0$).  Then the previous example says that this
completely characterizes the total wealth of the investor at any given
time, as long as the initial wealth is given.
\end{example}
A definition of state European and state American contingent claims is
given below.  It is known that state completeness is equivalent to the
maximality of the range of the matrix $\sigma(s_0,t,p^0)$ a.e. with
respect to the Lebesgue measure for all $s_0\leq t \leq T$.
See~\citet{Londono04}.
 
\section{State tameness and state arbitrage.  Characterization}\label{sec:statetameness}
\begin{definition}\label{D:w-tame}
  A $((\pi_0,\pi),\Gamma,\tau)$ portfolio evolution structure with
  random time $\tau$ after time $s_0$ financed by the income $\Gamma$
  is said to be \emph{state-tame}, if the process
  $H(s,\cdot,x,p)G(s,\cdot,x,p)$ is uniformly bounded below for all
  $x\in\mathbb{R}$, $p\in\mathbb{R}^{n+1}_+$, (the bound could depend
  on $x$, $p$ and $s$), where the process $G$ is defined by
\begin{multline}
  G(s,t,x,p)=B(s,t,p)\int_s^tB^{-1}(s,u,p)\pi^{\prime}(s,u,x,p)\sigma(s,u,p)\,dW_s(u) +\\
  B(s,t,p)\int_s^tB^{-1}(s,u,p)\pi^{\prime}(s,u,x,p)\left(b(s,u,p)+\delta(s,u,p)-r(s,u,p)\mathbf{1}_n\right)\,du
\end{multline}
\end{definition}
\begin{remark}
  We point out that equation \eqref{eq:wealthportfolio_process}
  implies that
\[
G(s,t,x,p)=X(s,t,x,p)-xB(s,t,p)-B(s,t,p)\int_s^tB^{-1}(s,u,p)\,d\Gamma(s,u,x,p)
\]  
In other words, $G$ is the process that measures at time $s$ the value
of the total gain using the portfolio $(\pi_0,\pi)$ after paying
interest at a rate given by the interest rate process $r$ for the
initial capital $x$, and after discounting the amount of money that a
bank account would have paid if the income stream would have been
saved in an interest rate account.  Therefore it is natural to call
the process $G$ the $\emph{gain in excess process}$.  In case $x=0$
and $\Gamma\equiv 0$, we obtain a gain in excess process for the self-
financed portfolio as described by \citet{Karatzas98}.
\end{remark}
\begin{definition}\label{def:Viability}
  A self-financed state-tame portfolio evolution structure
  $((\pi_0,\pi),\tau)$ with random time $\tau$ is said to be a
  \emph{state arbitrage opportunity for the initial wealth $x$ and
    initial price configuration $p\in\mathbb{R}_+^{n+1}$} if
\begin{equation*}
  \mathbf{P}\left[H(s,t,p)G(s,t,x,p)\geq 0\right]=1,\quad\text{and}\quad\mathbf{P}\left[H(s,t,p)G(s,t,x,p)>0\right]>0
\end{equation*}
where $G$ is the gain in excess process that corresponds to
$((\pi_0,\pi),\tau)$. We say that a market $\mathcal{M}$ is
\emph{state-arbitrage-free} if there are not portfolio evolution
structures with wealth $x$ and price configuration $p$ that are state
arbitrage opportunities.
\end{definition}

\begin{theorem}\label{T:Viability}
  A market $\mathcal{M}$ is state-arbitrage-free if and only if the
  $P$-consistent family of processes $\theta$ for all $s_0\leq s\leq
  T$, $p\in\mathbb{R}_+^{n+1}$ satisfies
\begin{equation}\label{E:Viability}
b(s,t,p)+\delta(s,t,p)-r(s,t,p)\mathbf{1}=\sigma(s,t,p)\theta(s,t,p)\qquad \mu_s\text{a.s.}
\end{equation}
\end{theorem}
\begin{remark}\label{R:localmartingale}
  We observe that if the family of processes $\theta$ satisfies
  equation$~(\ref{E:Viability})$ then for any initial capital $x$,
  initial price $p$, and wealth income evolution structure
  $(X,\Gamma,\tau)$,
\begin{eqnarray}\label{eq:localmartingale}
\lefteqn{Y(s,t,x,p)\equiv H(s,t,p)X(s,t,x,p)-\int_{s}^t H(s,t,p)\,d\Gamma(s,u,x,p) }\\
& =x+\int_s^tH(s,u,p)\left[\sigma^{\prime}(s,u,p)\pi(s,u,x,p)-X(s,u,x,p)\theta(s,u,p)\right] ^{\prime}\,dW_s(u). &\nonumber
\end{eqnarray}
As a straightforward consequence of It\^o's calculus, a portfolio
evolution structure $((\pi_0,\pi),\Gamma,\tau)$ is state tame if and
only if the process defined by the expression given on the left hand
side of the last equation is uniformly bounded below (where the bound
might depend on $x$, $p$, and $s$).  In case $Y$ is uniformly bounded
below, as it happens when $(X,\Gamma,\tau)$ is a hedgeable wealth
income structure with a state tame portfolio, $Y(s,\cdot,x,p)$ is a
super-martingale, for all $x$, $p$, and $s$.  The latter follows as
a result of $Y(s,\cdot,x,p)$ being a local-martingale that is uniformly
bounded below. Hence under these conditions,
\[
x\geq \mathbf{E}\left[Y(\tau)(x,p)\right]>-\infty.
\]
\end{remark}
\begin{proof}[Proof of Theorem \ref{T:Viability}]
  First, we prove necessity. For $s_0\leq s\leq t\leq T$, let us
  define the $P$-consistent progressive measurable process with two
  parameters of class
  $C^{0,\epsilon}(\mathbb{R}_+^{n+1}\colon\mathbb{R}_+^n)$ (where
  $\epsilon$ is an appropriate positive constant),
\[
\kappa(s,t,p)= b(s,t,p)+\delta(s,t,p)-r(s,t,p)\mathbf{1}_n -
\sigma(s,t,p)\theta(s,t,p),
\]
Define $\pi_0$ and $\pi$ to be the $P$-consistent processes of class
$C^{0,\epsilon}$
\begin{gather*}
  \pi(s,t,p)=\kappa(s,t,p)\notag\\
  \pi_0(s,t,p)=B(s,t,p)\int_s^t
  B^{-1}(s,u,p)\kappa^{\prime}(s,u,p)\kappa(s,u,p)\,du
  -\kappa^{\prime}(s,t,p)\mathbf{1}_n
\end{gather*}
It follows that $(\pi_0,\pi)$ is a self-financed portfolio evolution
structure with gain in excess process given by
\[
G(s,t,p)=\int_s^tB^{-1}(s,u,p)\kappa^{\prime}(s,u,p)\kappa(s,u,p)\,du
.
\]
Since $H(s,t,p)G(s,t,p)\geq 0$, the non-state-arbitrage hypothesis
implies the desired results.  To prove sufficiency, assume that the
family of processes $\theta$ satisfies equation \eqref{E:Viability}
for all $s$, $p$, and let $(\pi_0,\pi)$ be a self-financed portfolio
structure with gain in excess process $G$.  Remark
\ref{R:localmartingale} implies that $H[p](s,\cdot)G[x,p](s,\cdot)$ is
a local-martingale for all $x$, $p$, and $s$.  State-tameness implies
that it is also bounded below. Fatou's lemma implies that
$H(s,\cdot,p)G(s,\cdot,x,p)$ is a super-martingale.  The result follows.
\end{proof}

\section{A view of state European contingent claims.}\label{sec:europeanvaluation}
Throughout the rest of the paper we assume that
equation~(\ref{E:Viability}) is satisfied for all $s$, $p$,
$\mu_s$-almost surely.

We propose to extend the concepts of European contingent claim,
hedgeability and completeness within the framework proposed.
\begin{definition}\label{D:contingent_claim}
  A \emph{state European contingent claim (SECC) with expiration date
    $\tau$} is a wealth-income evolution structure
  $(X,\Gamma,\tau)\in\mathcal{X}(\mathcal{M})$ where the family of
  processes defined by
  \begin{equation}\label{eq:value_process}
Y(s,t,x,p)\equiv H(s,t,p)X(s,t,x,p)-\int_{s}^t
H(s,t,p)\,d\Gamma(s,u,x,p)    
\end{equation}  
are uniformly bounded from below continuous semimartingales for all
$x$, $p$ and $s$ (where the bound might depend on $x$, $p$, and $s$).
Moreover it is assumed that
\begin{equation}\label{E:valuesecc}
 x= \mathbf{E}\left[Y_s(\tau)(x,p)\right].
\end{equation}
We shall say that $Y$ is the \emph{discounted payoff process after
  time $s$ for the SECC}.
\end{definition}
\begin{remark}\label{rem:value_martingale}
  Under the conditions of Definition \ref{D:contingent_claim},
  equation \eqref{E:valuesecc} is equivalent to require that the
  process defined by equation \eqref{eq:value_process} is a martingale
  for each $x$, $p$, and $s$.
\end{remark}
\begin{definition}\label{D:wcompleteness}
  A state European contingent claim $(X,\Gamma,\tau)$ is called
  \emph{hedgeable} if $(X,\Gamma,\tau)$ is a hedgeable wealth-income
  evolution structure (by a state-tame portfolio evolution structure).
  The market model $\mathcal{M}$ is called \emph{state complete} if
  every state European contingent claim is hedgeable. Otherwise it is
  said to be \emph{state incomplete}.
\end{definition}

For the following theorem we assume $\left\{i_1<\cdots <
  i_k\right\}\subseteq \left\{1,\cdots, d\right\}$ is a set of indexes
and let $\left\{i_{k+1}<\cdots<i_{d}\right\}\subseteq
\left\{1,\cdots,d\right\}$ be its complement. Let $\sigma_i(s,t,p)$,
$1\leq i\leq k$, be the $i^{th}$ column process for the matrix valued
process $(\sigma_{i,j}(s,t,p)), s\leq t\leq T$.  Namely,
$\sigma_i(s,t,p)$, for $1\leq i\leq k$, is the $\mathbb{R}^n$-valued
progressively measurable process whose $j^{th}$ entry, for $1\leq
j\leq d$ agrees with $\sigma_{i,j}(s,t,p)$, for $0\leq s\leq t\leq
T$.  We denote by $\sigma_{i_1,\cdots,i_k}(s,t,p)$, $0\leq t\leq T$
the $n\times k$ matrix valued process whose $j^{th}$-column process
agrees with $\sigma_{i_j}(s,t,p)$, $0\leq s\leq t\leq T$, for $1\leq
j\leq k$.  We  adopt the same notation for other processes: if
$i_1<\cdots<i_{k^{\prime}}$ is a set of indexes, and $\psi_{i_j}$ for
$1\leq j\leq k^{\prime}$ is a process taking values in some euclidean
space $\mathbb{R}^{n^{\prime}}$, then we  denote by
$\psi_{i_1,\cdots,i_{k^{\prime}}}$ the $k^{\prime}\times n^{\prime}$
matrix valued process, where the $j^{th}$-column process agrees with
$\psi_{i_j}$.  We  denote by
$(\mathcal{F}_{s,t}^{i_1,\cdots,i_k})$ the filtration with two
parameters that is the $\mathbf{P}$ augmentation by the null sets of
the natural filtration
$\left\{\sigma(W^{i_1}_s(t),\cdots,W^{i_k}_s(t), 0\leq s\leq t\leq
  T\right\}$.  We  say that a wealth and income evolution
structure $(X,\Gamma,\tau)$ is a
$(\mathcal{F}_{s,t}^{i_1,\cdots,i_k})$ wealth-income evolution
structure if $X(s,\cdot,x,p)$, $\Gamma(s,\cdot,x,p)$ are
$\mathcal{F}^{i_1,\cdots,i_k}_{s,t}$ progressive measurable processes
and $\tau_s(x,p)$ is a stopping time relative to the filtration
$\mathcal{F}^{i_1,\cdots,i_k}_{s,t}$ for all $x$, $p$, and $s$. In
addition if $(X,\Gamma,\tau)$ is a SECC then we  say that it is a
\emph{$\mathcal{F}^{i_1,\cdots,i_k}_{s,t}$ SECC}.  First we prove a
lemma that would be needed for the proof of necessity in
Theorem~\ref{thm:completeness}.  Compare with \citet[Lemma
1.6.9]{Karatzas98}

\begin{lemma}\label{lemma:c-infinity}
  There exists a function $\varphi\colon
  L(\mathbb{R}^n\colon\mathbb{R}^k)\to\mathbb{R}^k$ of class
  $C^{\infty}$ when $L(\mathbb{R}^n\colon\mathbb{R}^k)$ is identified
  with $\mathbb{R}^{n\times k}$ in the natural way, with the property
  that $\varphi(\sigma)\in Im(\sigma^{\prime})$, $\varphi(\sigma)\neq
  \mathbf{0}$ if $Im(\sigma^{\prime})\neq\{\mathbf{0}\}$.
\end{lemma}
\begin{proof}
  Let $j_1<\cdots <j_r$ be a set of indexes in $\{1,\cdots,k\}$.
  Define $D_{j_1\cdots,j_r}\subset L(\mathbb{R}^n\colon\mathbb{R}^k)$
  to be a set of matrices $\sigma$ such that
  $\{\sigma_{j_1}^{\prime}\cdots,\sigma_{j_r}^{\prime}\}$ is linearly
  independent,
  $Im(\sigma_{j_1,\cdot,j_r}^{\prime})=Im(\sigma^{\prime})$, and there
  is not a set $j_1^{\prime}<\cdots<j^{\prime}_r$ with
  $j_{r}^{\prime}=j_r$ of indexes with
  $Im(\sigma^{\prime}_{j_1^{\prime},\cdots,j_r^{\prime}})=Im(\sigma_{j_1,\cdot,j_r}^{\prime})$.
  Define the function
  $\varphi^{j_1,\cdots,j_r}(\sigma)=\sum_{j_i}\omega_{i}(\sigma^{\prime}_{j_1}\cdots,\sigma^{\prime}_{j_i})\sigma^{\prime}_{j_i}$
  on $D_{j_1,\cdots,j_r}$, where $\omega_{i}$ is an $i$-form that
  computes a signed volume of the simplex formed by the the convex
  combination of the vectors in its arguments.  Define $\varphi\mid
  D_{j_1,\cdots,j_r}=\varphi^{j_1\cdots j_r}$, for each sequence
  $j_1<,\cdots,<j_r$ of indexes in $\{1,\cdots, k\}$.  It follows that
  $\varphi$ is a $C^{\infty}$ function, that satisfies the required
  conditions.
\end{proof}

In the following we assume that $X_1(t),\cdots,X_k(t)$ is a sequence
of continuous $\mathcal{F}^{i_1,\cdots,i_k}_{0,t}$ semimartingales
with values in $\mathfrak{X}^{0,\delta}(\mathbb{R}^n_+)$, the linear
topological space of $C^{0,\delta}$ vector fields for some $\delta>0$;
see \citet[Sec. 4.2]{hK90} for a description of semimartingales with
values in $\mathfrak{X}^k(M)$ for any non-negative integer $k$, and a
manifold $M$.  In addition, we also assume that
$X_1(t)-X_1(s),\cdots,X_k(t)-X_k(s)$ are
$\mathcal{F}_{s,t}^{i_1,\cdots,i_k}$ semimartingales.
\begin{theorem}\label{thm:completeness}
  Assume that for each $s$ and $p$, $\theta_i(s,t,p)=0$ $\mu_s$ a.s.
  for $i\notin\left\{ i_1,\cdots,i_k\right\}$ where
  $\theta(s,t,p)=(\theta_1(s,t,p),
  \cdots,\theta_d(s,t,p))^{\prime}$ is the market price of risk.
  Assume that $\sigma_{i_1,\cdots,i_k}(s,t,p)$ is a matrix valued
  process of two parameters such that
  \begin{multline*}
    Range^{\perp}(\sigma_{i_{k+1},\cdots,i_{d}}(s,t,p))=Range(\sigma_{i_1,\cdots,i_k}(s,t,p))\\
    =gen\{X_1(t)\cdots,X_k(t)\}
  \end{multline*}
  a.s. $\mu_s$ for all $p$, $s$. In addition assume that the interest
  rate process $B$ is a $(\mathcal{F}^{i_1,\cdots,
    i_k}_{s,t})$-process with two parameters.  Then, any
  $\mathcal{F}^{i_1,\cdots,i_k}_{s,t}$ state European contingent claim
  is hedgeable if and only if
  $Rank(\sigma_{i_1,\cdots,i_k}(s,t,p))=k$ a.s. $\mu_s$ for all $p$,
  $s$.  In particular, a financial market $\mathcal{M}$ is state
  complete if and only if for all $s$ and $p$, $\sigma(s,t,p)$ has
  maximal range a.s. $\mu_s$.
\end{theorem}
\begin{proof}[Proof of sufficiency]
  Let $(X,\Gamma,\tau)$ be a $\mathcal{F}^{i_1,\cdots,i_k}_{s,t}$
  SECC.  As a consequence of equation~\eqref{E:valuesecc}, the process
  defined by equation \eqref{eq:value_process} is a martingale for all
  $x$, $p$, and $s$.  By \citet[Exercise 3.2.10]{hK90} there exists
  $\epsilon>0$, and a $(\mathcal{F}_{s,t}^{i_1,\cdots,i_k})$-process
  of class
  $C^{0,\epsilon}(\mathbb{R}\times\mathbb{R}^{n}_+\colon\mathbb{R}^d)$,
  $\varphi(s,t,x,p)=\left(\varphi_1(s,t,x,p),\cdots,\varphi_d(s,t,x,p)\right)^{\prime}$,
  $t\in\left[s,\tau_s(x,p)\right]$, such that
\begin{multline*}
  H(s,t,p)X(s,t,x,p)-\int_{s}^tH(s,u,p)\,
  d\Gamma(s,u,x,p)\ \\
  =x + \int_s^t\varphi ^{\prime}(s,u,x,p)\,dW_s(u)
\end{multline*}
for all $x$, $p$, $s$, and $t$, where $\varphi_i(s,t,x,p)=0$ for
$i\notin\left\{i_1,\cdots,i_k\right\}$.  Define $\pi(s,t,x,p)$,
$t\in\left[s,\tau_s(x,p)\right]$, as the unique
$C^{0,\epsilon}(\mathbb{R}\times\mathbb{R}^n_+\colon\mathbb{R}^n)$
process with values in $ker^{\perp}(\sigma^{\prime}(s,t,p))$ such
that
\begin{equation*}
  \sigma ^{\prime}(s,t,p)\pi(s,t,x,p)=H^{-1}(s,t,p)\varphi(s,t,x,p)+ X(s,t,x,p)\theta(s,t,p) .
\end{equation*}
The existence and uniqueness of such a portfolio follows from the
hypotheses. (For instance, use the Gauss-Jordan algorithm to obtain a
solution of the system
$\sigma_{i_1,\cdots,i_k}^{\prime}\pi_{i_1,\cdots,i_k}=H\varphi_{i_1,\cdots,i_k}+X\theta_{i_1\cdots,i_k}$
and next take the projection onto
$ker^{\perp}(\sigma^{\prime}_{i_1,\cdots,i_k})=Range(\sigma_{i_1,\cdots,i_k})$.)
Define
$\pi_0(s,t,x,p)=X(s,t,x,p)-\pi^{\prime}(s,t,x,p)\mathbf{1}_n$.  It
follows using It\^o's formula that $(\pi_0,\pi)$ is a portfolio
process that finances the wealth $X$ with income $\Gamma$.
\end{proof}

\begin{proof}[Proof of necessity]  Let
  $\varphi: L(\mathbb{R}^k;\mathbb{R}^n)\mapsto\mathbb{R}^k$ be a
  bounded $C^{\infty}$ function defined by the Lemma
  \ref{lemma:c-infinity}. Let us define $\psi(s,t,p)$ to be the
  bounded, $\mathcal{F}^{i_1,\cdots,i_k}_{s,t}$-progressively
  measurable process
  $\psi_{i_1,\cdots,i_k}(s,t,p)=\varphi(\sigma_{i_{k+1},\cdots,i_d}(s,t,p))$
  and $\psi_j(s,t,p)=0$ for $j\notin \left\{i_1,\cdots,i_k\right\}$.
  We consider the $\mathcal{F}^{i_1,\cdots,i_k}_{s,t}$-progressively
  measurable SECC with no income, expiration date $T$, and whose
  wealth process is defined by
\begin{equation*}
        X(s,t,x,p)=x+\int_s^t\frac{1}{H(s,u,p)}\psi ^{\prime}(s,u,p)\,
        dW_s(u),\qquad\mbox{for }s\leq t\leq T  .
\end{equation*}
Let $(\pi_0,\pi)$ be the state tame portfolio that finances the given
wealth. It follows that
\begin{equation*}
H(s,t,p)X(s,t,x,p)=x +\int_s^t\psi^{\prime}(s,u,p)\, dW_s(u)
\end{equation*}
is a martingale for all $x$, $p$ and $s$.  Using
equation~\eqref{eq:localmartingale}, and \citet[Exercise 3.2.10]{hK90}
we obtain
\begin{eqnarray*}
\psi_{i_1,\cdots,i_k}(s,t,p)=\sigma_{i_1,\cdots,i_k}^{\prime}(s,t,p)\pi_{i_1,\cdots,i_k}(s,t,x,p)-X(s,t,x,p)\theta_{i_1,\cdots,i_k}(s,t,p)   \nonumber\\
\in Ker^{\perp}(\sigma_{i_1,\cdots,i_k}(s,t,p)\cap Ker(\sigma_{i_1,\cdots,i_k}(s,t,p))=\left\{\mathbf{0}\right\} 
\end{eqnarray*}
a.s $\mu_s$ for all $s$. The result follows.
\end{proof}

\section[A view of state American contingent claims.]{A view of state American contingent claims.}\label{sec:review-state-american}
\begin{definition}\label{def:sacc}
  Let $s_0\in\left[0,T\right]$. Assume a wealth-income evolution
  structure after time $s_0$,
  $(X,\Gamma,\tau)\in\mathcal{X}(\mathcal{M})$, such that the family
  of processes defined by equation~\eqref{eq:value_process} are
  uniformly bounded from below continuous semimartingales for all $x$,
  $p$ and $s$ (where the bound might depend on $x$, $p$, and $s$). A
  \emph{State American Contingent Claim} is a wealth income structure
  as above with the property that for all $x$, $p$ and $s$,
\begin{equation}\label{eq:valuesacc}
 x =\sup_{\tau^{\prime}\in\mathcal{S}_s(X,\tau)} \mathbf{E}\left[Y_s(\tau^{\prime})(x,p)\right] 
\end{equation}
where $\mathcal{S}_{s}(X,\tau)$ denotes, for a given wealth evolution
structure $(X,\tau)$ and a given $s$, the family of stopping times
that takes values in $[s,\tau\vee s]$ that are $(X,P)$-consistent
after time $s$ (here $s\vee t=\max(s,t)$ and it is assumed that
$\sup\emptyset =\infty$).  We  call the process $Y(s,\cdot,x,p)$
the \emph{discounted payoff process after time $s$}.  We  denote
the quantity on the right hand side of the equation
\eqref{eq:valuesacc} as $u_{s}(x,p)$ and we  say that $u_s(x,p)$
is the \emph{value of the consistent state American contingent claim
  at time $s$} (given that the price of the stock process at time $s$
is $p$ and the wealth process is $x$); in case $s=0$ the process
$Y(x,p)$, and $u_0(x,p)$ are simply called \emph{the discounted payoff
  process} and \emph{the value of the consistent state American
  contingent claim}, respectively.
  \end{definition}
  \begin{remark}\label{rem:supermartingale}
    Under the conditions of Definition \ref{def:sacc},
    $(X,\Gamma,\tau)\in\mathcal{X}(\mathcal{M})$ satisfies
    equation~\eqref{eq:valuesacc} iff the family of processes
    $Y(s,\cdot\wedge\tau,x,p)$ are super-martingales for all $x$, $p$
    and $s$.
  \end{remark}

  In addition to properties mentioned above we  assume for the
  remainder of this section the following condition.
  \begin{condition}\label{con:continuity}
    For any $s_0\leq s\leq T$ and all stopping times
    $\tau^{\prime}\in\mathcal{S}_s(X,\tau)$, the function
\begin{equation}
  \label{eq:continuity}
\varphi_{s,\tau^{\prime}}(x,p)=\mathbf{E}\left[Y_s(\tau^{\prime}\wedge \tau)(x,p)\right] 
\end{equation}
is a continuous function in $(x,p)$, and the given family of functions
is a equicontinuous set of functions on compact sets (in $(x,p)$).
Moreover assume that there exist positive constants $\gamma\geq 1$,
$\alpha_1,\alpha_2,\alpha_3,\beta_0,\cdots,\beta_n$, with
$\alpha_1^{-1}+\alpha_2^{-1}+\alpha_3^{-1}+\sum_{i=0}^n\beta_i^{-1}<1$ such
that the random field $Y(s,t,x,p)\triangleq Y_s(\tau)(s,t,x,p)$
satisfies
\begin{multline}
  \label{eq:continuity_random}
  \mathbf{E}\left[\mid
    Y(x,p,s,t)-Y(x^{\prime},p^{\prime},s^{\prime},t^{\prime})\mid^{\gamma} \right]\leq \\
  C\left(\mid s-s^{\prime}\mid^{\alpha_1}+\mid
    t-t^{\prime}\mid^{\alpha_2}+ \mid x-x^{\prime}\mid^{\alpha_3}
    +\sum_{i=0}^n\mid p_i-p_i^{\prime}\mid^{\beta_i}\right).
\end{multline}
  \end{condition}
  The previous condition is usually satisfied when $X$ is a process
  that solves a stochastic differential equation.  For instance, see
  Example \ref{ex:construction_wealth}, and \citet[Lemma
  4.5.6]{hK90}.  Equation (\ref{eq:continuity_random}) above is needed
  in order to obtain a continuous modification of the random field and
  its conditional expectation.  See Kolmogorov's continuity criterion
  of random fields (\citet[Theorem 1.4.1 and Exercise 1.4.12]{hK90}).
  In the following, conditional expectations of stochastic processes
  are the continuous modifications of the given stochastic processes. 
  
  In order to state and prove a theorem for valuation of State
  American contingent claims is needed to define when a wealth income
  evolution structure outperforms other.  
  \begin{definition}
    We  say that a wealth and income evolution structure
  $(X^{\prime},\Gamma^{\prime}, \tau^{\prime})$ dominates a wealth income
  structure $(X,\Gamma,\tau)$ if $\tau^{\prime}_s(x,p)\geq
  \tau_s(x,p)$ for all $s$, $x$ and $p$, 
  \begin{eqnarray*}
X^{\prime}(s,t,x,p)\geq X(s,t,x,p)\qquad\text{and}\qquad 
\Gamma^{\prime}(s,t,x,p)\leq \Gamma(s,t,x,p)    
  \end{eqnarray*}
for all $s\leq t\leq \tau_s(x,p)$.
  \end{definition}
  
Theorem \ref{thm:american} below provides conditions under which 
  every state American contingent claim is dominated by a hedgeable
  state American contingent claim; compare with Theorem
  5.1 \citet{Londono04}.
\begin{theorem}\label{thm:american}
  Assume that the hypotheses of Theorem \ref{thm:completeness} and
  Condition \ref{con:continuity} hold.  Then, any
  $\mathcal{F}^{i_1\cdots,i_k}_{s,t}$-progressively measurable state
  American contingent claim $(X,\Gamma,\tau)$ is dominated by a hedgeable $\mathcal{F}^{i_1\cdots,i_k}_{s,t}$-progressively measurable state
  American contingent claim  if and only
  if for each $s$ and $p$, $Rank(\sigma_{i_1,\cdots,i_k}(s,t,p))=k$
  a.s. $\mu_s$.  In particular, a financial market $\mathcal{M}$ is
  American state complete if and only if $\sigma(s,t,p)$ has maximal
  range a.s.  $\mu_s$ for all $s$.
\end{theorem}
The sufficiency is a consequence of Remark \ref{rem:supermartingale},
and Theorem \ref{thm:completeness}; the proof of necessity is given
below.

First we point out an elementary fact, needed for the proof.  If
$\sigma\in\mathcal{S}_{t}(X,\tau)$ for $t\geq s$, then $\sigma$ can be
seen as an element of $\mathcal{S}_s(X,\tau)$ using the natural
identification, namely $\sigma$ is identified with the stopping
structure $\sigma^{s}\in\mathcal{S}_s(X,\tau)$ defined by
$\sigma_{t^{\prime}}^s(x,p)=\sigma_{t^{\prime}}(x,p)$ for
$t^{\prime}\geq t$, and $\sigma_{t^{\prime}}(x,p)=\sigma_{t}(X(t^{\prime},t,x,p),P(t^{\prime},t,x,p)]$
otherwise.  We  denote both elements in the same way, hoping that the meaning
will be clear from context.
\begin{lemma}\label{lem:american_aproximation}
  Assume Condition \ref{con:continuity} and let  $(X,\Gamma,\tau)$ be a state American contingent claim and let
  $s_0\leq s\leq T$.  Assume $\tau_1,\tau_2
  \in\mathcal{S}_{s}(X,\tau)$.  Then, there exists
  $\tau^{\prime}\in\mathcal{S}_s(X,\tau)$ with the property that for
  any $s\leq s^{\prime}\leq T$,
\begin{equation*}
u_{s^{\prime}}(x,p)\geq\mathbf{E}[Y_{s^{\prime}}(\tau^{\prime})(x,p)]\geq
\max\{\mathbf{E}[Y_{s^{\prime}}(\tau_1)(x,p)],\mathbf{E}[Y_{s^{\prime}}(\tau_2)(x,p)]
\end{equation*}
and
\begin{equation*}
  \mathbf{E}\left[Y_{s^{\prime}}\left(\tau^{\prime}\right)(x,p)\mid
    \mathcal{F}_{s^{\prime},t}\right]\geq
  \max\left\{\mathbf{E}\left[Y_{s^{\prime}}(\tau_1)(x,p)\mid
      \mathcal{F}_{s^{\prime},t}\right],
    \mathbf{E}\left[Y_{s^{\prime}}\left(\tau_2\right)(x,p)\mid\mathcal{F}_{s^{\prime},t}\right]\right\}.
\end{equation*}
a.s. $\mu_{s^{\prime}}$.
\end{lemma}
\begin{proof}
  Define
\begin{eqnarray*}
    \label{eq:1}
\lefteqn{\tau_{s^{\prime}}(x,p)=(\tau_1\wedge\tau_2)_{s^{\prime}}(x,p)\mathbf{1}_{\mathbf{E}\left[Y_{s^{\prime}}\left(\tau_1\vee\tau_2\right)(x,p)\mid\mathcal{F}_{s^{\prime},t}\right]\left(\tau_1\wedge\tau_2\right)_{s^{\prime}}(x,p)<Y_{s^{\prime}}\left(\tau_1\wedge\tau_2\right)(x,p)}}\\
& &\mbox{} + (\tau_1\vee\tau_2)_{s^{\prime}}(x,p)\mathbf{1}_{\mathbf{E}\left[Y_{s^{\prime}}\left(\tau_1\vee\tau_2\right)(x,p)\mid\mathcal{F}_{s^{\prime},t}\right]\left(\tau_1\wedge\tau_2\right)_{s^{\prime}}(x,p)\geq
  Y\left(\tau_1\wedge\tau_2\right)_{s^{\prime}}(x,p)} 
  \end{eqnarray*}
  for all $s^{\prime}\geq s$ where $t\vee
  t^{\prime}=\max\left(t,t^{\prime}\right)$, and $t\wedge
  t^{\prime}=\min\left(t,t^{\prime}\right)$.  Then $\tau$ has the required properties.
\end{proof}

\begin{proof}[Proof of necessity.]
  Let $Y$ be the discounted payoff structure of a state American
  contingent claim $(X,\Gamma,\tau)$.  For each $x$, $p$ and $s$ there
  exists a sequence of families of consistent stopping times
  $\sigma_n^{x,p}\in\mathcal{S}_{s}(X,P)$ along sets $N_s\in
  \mathcal{P}_s$ of zero $\mu_s$ measure such that $
  \mathbf{E}\left[Y_s\left(\sigma_n^{x,p}\right)(x,p)\right] \uparrow
  x=u_s(x,p)$, with
\[
\mathbf{E}\left[Y\left(\sigma_{n+1}^{x,p}\right)(x,p)\mid\mathcal{F}_{s,t}\right]\geq\mathbf{E}\left[Y\left(\sigma_n^{x,p}\right)(x,p)\mid\mathcal{F}_{s,t}\right]\qquad\text{
  for }(t,\omega)\notin N_s.
\]  
The latter follows by Lemma \ref{lem:american_aproximation}.  Without
loss of generality the version chosen for the given random fields are
continuous (see \citet[Exercise 1.4.12]{hK90}).  Using a
diagonalizing process and Lemma \ref{lem:american_aproximation} again,
it is possible to prove that there exists a sequence of families of
consistent stopping times $\sigma_n\in\mathcal{S}_{s_0}(X,P)$ and a
set with the property that for any triple of points with rational
coordinates $x\in\mathbb{Q}$, $p\in\mathbb{R}_+^{n}\cap\mathbb{Q}^n$,
and $s\in [s_0,T]\cap\mathbb{Q}$ there exists $M$ (depending on $x$,
$p$, and $s$) large enough with
\[
\mathbf{E}\left[Y_s\left(\sigma_{n+1}\right)(x,p)\mid\mathcal{F}_{s,t}\right]\geq\mathbf{E}\left[Y_s\left(\sigma_{n}\right)(x,p)\mid\mathcal{F}_{s,t}\right],\qquad\text{ for
}(t,\omega)\notin N_{s}
\]  
where $\mu_s(N_s)=0$, and
\[
\mathbf{E}\left[Y_s(\sigma_n)(x,p)\right] \uparrow x=u_s(x,p).
\] 
Using Condition \ref{con:continuity}, and the Ascoli-Arzelà Theorem,
there exist a subsequence of  stopping times (that we also denote as $\sigma_n$) with $\sigma_{\tau_n}\to u$ uniformly on compact sets,
where $u(s,x,p)\equiv u_s(x,p)$ is the value of the SACC, and
$\varphi_{\sigma_n}(s,x,p)\equiv\varphi_{s,\sigma_n}(x,p)$.  Doob's
inequality shows that the stochastic process
$(\mathbf{E}\left[Y(\sigma_n)(x,p)\mid\mathcal{F}_{s,t}\right])_{s\leq
  t\leq T}$ is a Cauchy sequence in the sense of uniform convergence in probability uniformly
on compact sets.  By completeness of the space of local-martingales,
for each $s$ there exists a local-martingale $\overline{Y}(s,t,x,p)$,
$t\in\left[s,T\right]$, such that $\mathbf{E}\left[Y(\sigma_n)(x,p)\mid\mathcal{F}_{s,t}\right]\rightarrow\overline{Y}(s,t,x,p)$,
$t\in\left[s,T\right]$, uniformly in probability, (and the convergence
is uniform in compact sets).  It follows by continuity that
$\overline{Y}(s,t,x,p)\geq Y(s,t,x,p)$ for all $s$, $x$, $p$, a.s.
$\mu_s$, and clearly $\overline{Y}(x,p)(s,s)=u_s(x,p)$.  Define
$\tau_n$ to be the first hitting time of $\overline{Y}(s,t,x,p)$,
$t\in\left[s,\tau\right]$, to the set $\left[-n,n\right]^c$. 
Using \citet[Exercise 3.2.10]{hK90} it follows that there exists a 
$\mathbb{R}^d$-valued process of class $C^{0,\varepsilon}$ for some
$\varepsilon>0$, 
$\varphi(s,t,x,p)=\left(\varphi_1(s,t,x,p),\cdots,\varphi_d(s,t,x,p)\right)^{\prime}$,
$t\in\left[s,\tau_n\right]$, such that
\begin{equation*}
  \overline{Y}(s,t,x,p)=x + \int_{s}^t\varphi ^{\prime}(s,t,x,p)\,dW_s(u)
\end{equation*}
where $\varphi_i(s,t,x,p)=0$ for
$i\notin\left\{i_1,\cdots,i_k\right\}$.  Define $\overline{X}(s,t,x,p)$ by
\[
  H(s,t,p)\overline{X}(s,t,x,p)-\int_{s}^tH(s,u,p)\,
  d\Gamma(s,u,x,p)=\overline{Y}(s,t,x,p)
\]
and  $\pi(s,t,x,p)$, $s\leq
t\leq \tau$, as the unique $\mathbb{R}^n$-continuous  process such that
\begin{equation*}
  \sigma^{\prime}(s,t,p)\pi(s,t,x,p)=H^{-1}(s,t,p)\varphi(s,t,x,p)+ \overline{X}(s,t,x,p)\theta(s,t,p) .
\end{equation*}
The existence and uniqueness of such a portfolio follows from the
hypotheses (see Lemma 1.4.7 in \citet{Karatzas98}).  Define
$\pi_0(s,t,x,p)=\overline{X}(s,t,x,p)- \pi^{\prime}(s,t,x,p)\mathbf{1}_n$.
Using It\^o's formula it follows that $(\overline{X},\Gamma,\tau)$ is
a wealth income process with the desired characteristics.
\end{proof}  
\begin{remark}
  We point out that the definition of consistent state American
  contingent claim and state American contingent claim are not
  equivalent.  In fact the supremo that defines the value of state
  American contingent claim is greater or equal than the value of the
  consistent American contingent claim defined by
  \eqref{eq:valuesacc}, and \emph{a priori} could be strictly greater.  The
  latter follows because the class of stopping times contains the
  class of $(X,\Gamma,\tau)$-consistent families of stopping times.
\end{remark}

\bibliographystyle{plainnat}
\bibliography{utility,eafit,tameness}
\end{document}